\documentclass[12pt]{article}
\usepackage[cp1251]{inputenc}
\usepackage[T2A]{fontenc}
\usepackage[english]{babel}
\usepackage{amsfonts,amsmath}
\usepackage{amsmath,amssymb}
\usepackage{graphicx,fancybox}
\usepackage{epstopdf}
\usepackage{enumitem}
\usepackage[left=1.5cm, right=1.5cm, top=2cm, bottom=2cm, headsep=0.7cm, footskip=1cm
]{geometry}
\newtheorem{theorem}{Theorem}
\newtheorem{corollary}{Corollary}
\newtheorem{remark}{Remark}

\newtheorem{example}{Example}

\begin{document}
	
	\title {\hbox{\normalsize }\hbox{}Integrable Billiards and Related Topics}
	
	 \author{Misha Bialy, \ Andrey E. Mironov}	
	\date{}
	\maketitle
	
	\sloppy
	\begin{abstract}
This paper surveys our results on integrable billiards. We consider various models of billiards, including Birkhoff, outer, magnetic, and Minkowski billiards. Also, we discuss wire billiards and billiards in cones. For four models of convex plane billiards, we also discuss an isoperimetric-type inequality for the Mather $\beta$-function. We conclude with a section of open questions on this subject.
	
	\end{abstract}

\section{Introduction.}
In this survey, we discuss the question of integrability of various billiard models and the related results. We refer to \cite{tabachnikov1995billiards, kozlov-treschev} for general introduction to the subject. 

The most classical model of mathematical billiards was introduced and studied by G. D. Birkhoff. Birkhoff and Poritsky conjectured that the only integrable Birkhoff billiards in the plane are ellipses. It was explicitly formulated in \cite{poritsky} and since then it remains unsolved in full generality.

We discuss various approaches to this conjecture in the subsequent sections. Our approach to the Birkhoff--Poritsky conjecture in this paper is related to the notion of {\it total integrability}, also called sometimes $C^0$-integrability.

Let $T$ be an area-preserving the twist map of the cylinder, and $U$ be an open subset of the cylinder. We say that $T$ is \textit{totally integrable in $U$} if $U$ admits a foliation by rotational (i.e., winding once around the cylinder) invariant curves of $T$.

This notion is related to E. Hopf's rigidity theorem, stating that a Riemannian 2-torus with no conjugate points is necessarily flat.

It became apparent that the ideas of Hopf's rigidity result can be extended to various types of billiards using their variational principles.
To that end, we need new tools and ideas. In particular, we use a non-standard generating function for billiards and carefully chosen weights in the integral-geometric part of the method.

We discuss Birkhoff billiards in greater detail in Section \ref{sect:birkhoff}.

We find that the Birkhoff--Poritsky conjecture is meaningful and interesting for four natural models of convex billiards in the plane. These four models are important for other geometric problems, e.g., for Length/Area approximations of convex bodies by inscribed/circumscribed extremal polygons, see \cite{toth,albers-tabachnikov,baranzini-bialy-sorrentino}. We discuss these models in Section \ref{sect:other}.

Moreover, we formulate the effective form of the rigidity results for Birkhoff billiards in Section \ref{sect:effective}.

In Section \ref{sect:symmetry}, we describe the rigidity phenomenon for billiards with a symmetry property.

In Section \ref{sect:beta}, we discuss the relation of the Mather $\beta$-function to the question of integrability. We also state an isoperimetric-type inequality for the Mather $\beta$-function for billiards which was recently discovered in  \cite{baranzini-bialy-sorrentino}.

Results on magnetic and Minkowski billiards are presented in Section \ref{sect:magnetic}.

Very recent results on Birkhoff billiards in cones and wire billiards are discussed in Section \ref{sect:wire}.

We conclude with open questions in Section \ref{sect:questions}.


\subsection*{Acknowledgments.} It is a pleasure to thank Sergei Tabachnikov for the illuminating discussions we had for many years on mathematical billiards. We are also grateful to Konstantin Khanin for useful consultations on various aspects of dynamical systems. 

M. Bialy was partially supported by ISF grant 974/24, A. E. Mironov was supported by the state contract of the Sobolev Institute of Mathematics, project no. FWNF-2026-0026.

\section{Birkhoff billiards.} \label {sect:birkhoff}It is crucial to specify what is understood by integrability. Accordingly, there are various approaches to the Birkhoff--Poritsky conjecture.

A very beautiful approach, connected to algebraic geometry, assumes that there exists a first integral of the billiard ball motion which is polynomial in the velocities. This class of integrals is very natural from the point of view of classical mechanics, where all conserved quantities appear to be polynomial in momenta variables.
This algebraic approach to the Birkhoff--Poritsky conjecture was initiated by S.~Bolotin \cite{bolotin,bolotin2} and then studied in our papers \cite{BM1,BM2}.
Finally, A. Glutsyuk \cite{g1,g2}, using the results of \cite{BM1,BM2}, completed the proof of the algebraic Birkhoff--Poritsky conjecture for billiards on the plane and constant curvature surfaces.
Let us remark that this approach does not require closeness to the ellipses, although it is restricted to the class of algebraic curves and polynomial integrals. Moreover, this approach allows to consider piece-wise smooth boundaries.

Other important developments in the direction of the Birkhoff--Poritsky conjecture are the following.
In \cite {delshams}, it is shown that perturbations of ellipses create a splitting of separatrices.

In \cite{Innami} (see also \cite{BA} for a geometric approach), it is proved that if there exists a sequence of convex caustics with the rotation numbers tending to $1/2$, then the billiard table is an ellipse. In \cite{treschev},  evidence of possible integrable dynamics around a 2-periodic orbit is given. In \cite{Marco}, a polynomial entropy approach to the problem is suggested.

Finally, let us cite here a series of influential results by V. Kaloshin et al. (see e.g. \cite{kaloshin-huang-sorrentino, KS} and also a very recent \cite{koval}) proving a local version of the Birkhoff--Poritsky conjecture in a neighborhood of ellipses in  suitable functional spaces.
We refer to an excellent recent survey \cite{Kaloshin-lecturenotes} with references therein, and also to
\cite{KKK} for a conjunction of the local and total integrability approaches.

Our approach to the Birkhoff--Poritsky conjecture is related to the notion of {\it total integrability} (a term suggested by A. Knauf). This notion is intimately related to E. Hopf's theorem claiming that a Riemannian 2-torus with no conjugate points is necessarily flat. As a corollary, the unit tangent bundle of a Riemannian 2-torus is foliated by incompressible invariant tori if and only if the Riemannian metric is flat.
In order to get the corollary, one needs to assert that all the lifts of the geodesics lying on such an invariant torus are length-minimizing and thus have no conjugate points. It is unknown if other integrable Riemannian metrics on the 2-torus, e.g., rotationally invariant metrics or Liouville metrics, can be distinguished in variational terms.

\begin{figure}[h]
	\centering
	\includegraphics[width=0.55\textwidth]{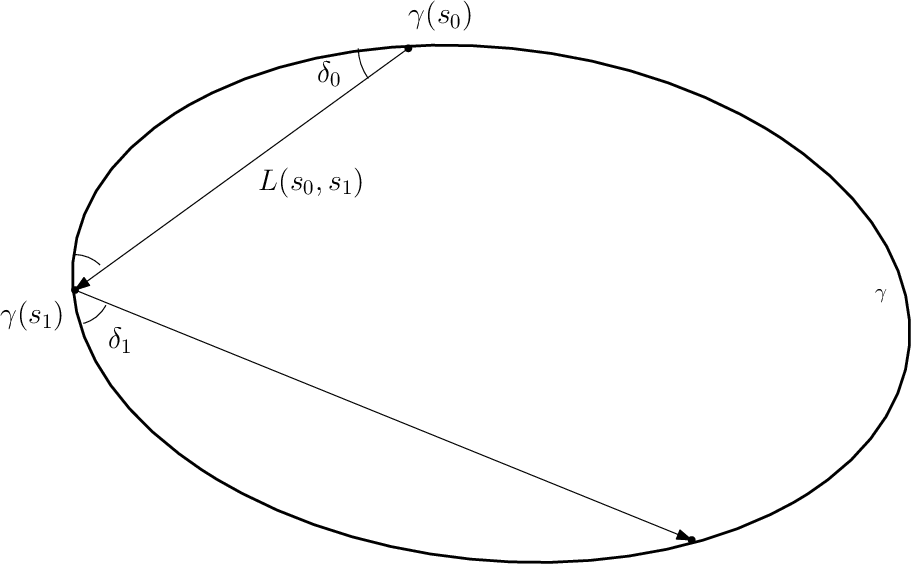}
	\caption{Generating function $L(s_0,s_1)=|\gamma(s_0)-\gamma(s_1)|$ }
	\label{fig:L}
\end{figure}

After analyzing E. Hopf's proof on geodesic flows with no conjugate points, it was understood in \cite {B0} (later in \cite{W} by a mirror formula) and then in \cite{B1} that Hopf's result can be generalized to convex billiards. This generalization is based on the fact that billiard orbits $\{s_n, n\in\mathbb Z\}$ are stationary sequences of the {\it Length} functional (see Fig. \ref{fig:L}).
\begin{equation}\label{functional}
	\mathcal L\{s_n\}=\sum_{n} L(s_n,s_{n+1}),\ L(s_n,s_{n+1})=|\gamma(s_n)-\gamma(s_{n+1})|,
\end{equation}

As a corollary, one gets that the only {\it totally integrable} billiards are circular.
Thus, the circular billiards are analogous to the geodesic flow on a flat torus.

\begin{figure}[h]
	\centering
	\includegraphics[width=0.55\textwidth]{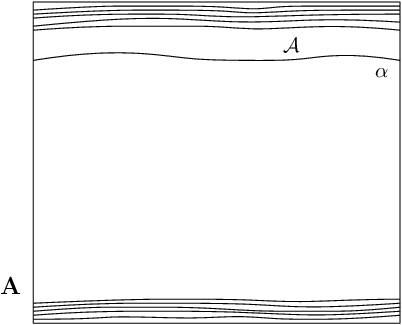}
	\caption{Invariant curve $\alpha$ of 4-periodic orbits and the region $\mathcal A$.}
	\label{fig:main-theorem}
\end{figure}

Surprisingly \cite{bialy-mironov}, elliptic billiards can be also distinguished among all centrally-symmetric billiards as  rigid objects in both: the variational and total integrability settings. The difference with circular billiards appears in the requirement of total integrability on a  certain neighborhood of the boundary of the phase cylinder $\mathbb{A}$, see Fig. \ref{fig:main-theorem}.
\medskip

\begin{theorem}\label{thm:birkhoff-1/4}
	Let $\gamma$ be a $C^2$-smooth, centrally-symmetric, simple closed curve of positive curvature. Suppose that the billiard ball map $T$ of $\gamma$  has a continuous rotational (i.e., winding once around the cylinder and simple) invariant curve $\alpha\in\mathbb A$ consisting of $4$-periodic orbits.
	
	Let $\mathcal A\subset \mathbb A$ be the domain of the phase cylinder between $\alpha$ and the boundary. Then the following two statements hold true:
	\begin{enumerate}
		\item	[(1)]If $\mathcal A$ is foliated by continuous rotational invariant curves, then $\gamma$ is an ellipse.
		
		\item	[(2)]If all orbits in $\mathcal A$ are locally maximizing for the Length functional, then $\gamma$ is an ellipse.
	\end{enumerate}
\end{theorem}
\medskip

Let us remark that the most general form of the Birkhoff--Poritsky conjecture assumes total integrability on any tiny neighborhood of the boundary. In Theorem \ref{thm:birkhoff-1/4}, the neighborhood is small but finite. Notice also that one cannot recover integrable billiards by infinitesimal data at the boundary. This fact follows from Marvizi--Melrose \cite{marvizi-melrose} results on interpolating hamiltonians.

Theorem \ref{thm:birkhoff-1/4} is restricted to the case of centrally-symmetric billiard tables. But the smoothness assumptions on the invariant curves are minimal (continuity and even less). Also, the billiard table is not assumed to be close to an ellipse.

Let us indicate the main novel tools of this approach. The first ingredient is a non-standard generating function for the billiard ball map leading to the twist property with respect to another vertical foliation of the phase cylinder $\mathbb A$.

The second ingredient is the remarkable structure of the invariant curve consisting of 4-periodic orbits.

The third ingredient is the use of integral geometry essentially as for the circular billiards, but now using the non-standard generating function, and carefully chosen weights. This is a new element in the integral geometry tool.

\section{Other convex plane billiards.}\label{sect:other}

We will summarize here four models of convex planar billiards.
All four models are very natural for other geometric problems, e.g., for Length/Area approximations of convex bodies by inscribed/circumscribed extremal polygons, see \cite{toth,albers-tabachnikov,baranzini-bialy-sorrentino} for further details.
These four models of billiards differ by being either inner or outer, and by a generating function of length or area.

Thus, classical Birkhoff billiards correspond to the choice of \emph{inner length} billiard.

The next famous example is the {\it outer billiard} (corresponding to the outer area billiard). Outer billiards were introduced by B. Neumann in 1959, but an earlier construction is due to M. Day in \cite{day}. Jurgen Moser popularized this system in the 1970's as a toy model for celestial mechanics: the orbit of a point around the  billiard table resembles the orbit of a celestial body around a planet or a star (see \cite{Moser, Moserbook}).

\begin{figure}[h]\label{fig:outer}
	\centering
	\includegraphics[width=0.65\linewidth]{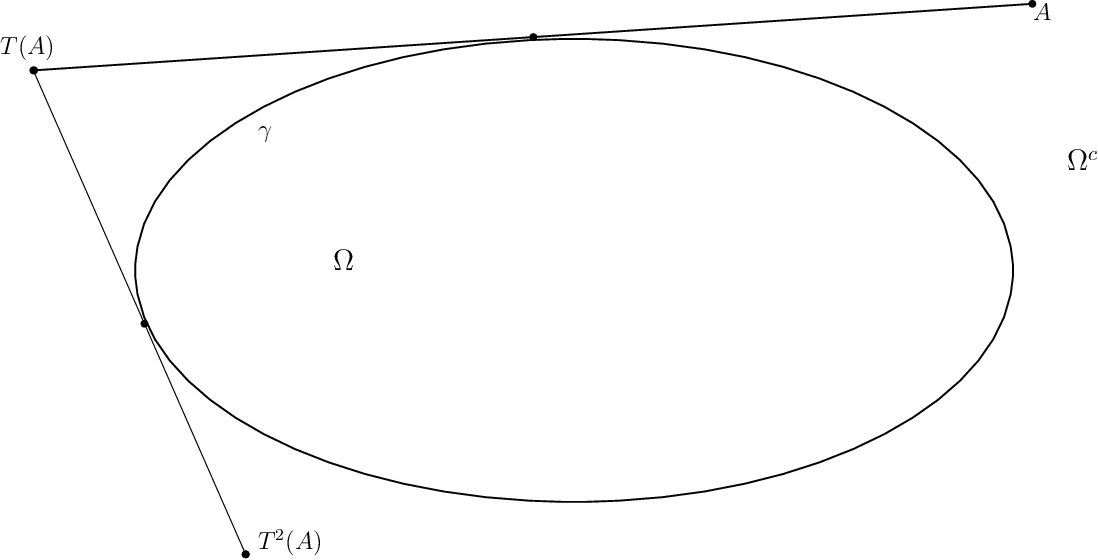}
	\caption{Outer billiard map}
\end{figure}

Consider a strictly convex domain $\Omega$ bounded by a closed curve $\gamma$ in the plane. Let $\Omega^c$ be the exterior of the domain. The outer billiard map $T$ acts on $\Omega^c$ as follows:

Given a point $A\in\Omega^c$, its image $T(A)$ is defined by the condition that the segment $[A,T(A)]$ is tangent to $\gamma$ exactly at the middle of the segment, Fig. \ref{fig:outer}. The map $T$ is a symplectic diffeomorphism
of $\Omega^c$ with respect to the standard form of the plane. Thus, $\Omega^c$ is the phase space of the outer billiard.

The following is a natural question which is analogous to the Birkhoff--Poritsky conjecture for usual billiards.
{\it Are there integrable outer billiards in the plane other than ellipses?}
For the algebraic version of this question, the answer is negative, as shown in  \cite{glutsyuk-shustin,tabachnikov-outer}.

It is easy to see from the definition of the outer billiard, that the group of affine transformations commutes with the outer billiard map. Thus, the outer billiard of the ellipse leaves an invariant foliation of $\Omega^c$ by concentric homothetic ellipses (since this is obviously the case for a circle).

For the total integrability of the outer billiard, the rigidity result remained resistant for a long time due to two main difficulties. First, the phase space $\Omega^c$ is not compact, which complicates the integral geometric part of the proof. Second, the affine nature of the problem makes it much harder than the euclidean one.

The following result as well as its variational version was proved in \cite{bialy-outer}.
\medskip

\begin{theorem}\label{thm:total}Let $\gamma$ be a $C^2$-smooth convex closed curve of positive curvature.
	Assume that the outer billiard of $\gamma $ is totally integrable, i.e., the phase space $\Omega^c$ is foliated by continuous rotational (i.e., non-contractable in $\Omega^c$) invariant curves, then $\gamma$ is an ellipse.
\end{theorem}
\medskip

Similarly to Theorem \ref{thm:birkhoff-1/4}, the proof of this result relies on two new ingredients: a \emph{non-standard} generating function of outer billiards, and the use of the suitable weights. Using them, the proof can be reduced to the Blaschke--Santalo inequality of affine geometry.

It is unknown if the rigidity result for the centrally-symmetric case, analogous to Theorem \ref{thm:birkhoff-1/4}, is valid.
Plausibly, the rigidity holds assuming: (a) the existence of an invariant curve $\alpha$ of $4$-periodic orbits, and (b) the total integrability of the outer billiard on the part of the phase space bounded by the invariant curve $\alpha$ and the boundary $\gamma$ (or for the part of the phase space which is the exterior of $\alpha$).

Notice that for a $C^2$-smooth curve $\gamma$, the existence of an invariant curve $\alpha$ consisting of $4$-periodic orbits is equivalent to $\gamma$ being a Radon curve \cite{Radoncurves} (see \cite{martini2006antinorms} and \cite{bialy2022self} for modern aspects).

\begin{figure}[h]
	\centering
	\includegraphics[width=0.5\linewidth]{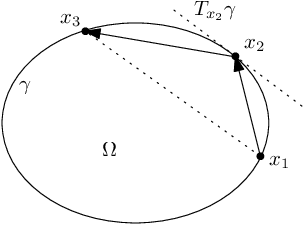}
	\caption{A symplectic billiard bounce $(x_1,x_2)\mapsto (x_2,x_3)$ in a strictly convex domain $\Omega$.}
	\label{figsimplbill}
\end{figure}

The next class with respect to the inner/outer and length/area classification is the class of \emph{symplectic billiards} (which are inner area billiards) introduced and studied in \cite{albers2018introducing}. This billiard bounces on the boundary $\gamma=\partial\Omega$ according to the following reflection law: given three successive impact points $x_1$, $x_2$ and $x_3$, the line joining $x_1x_3$ and the tangent line $T_{x_2}\gamma$ of $\partial\Omega$ at $x_2$ are parallel (see Fig. \ref{figsimplbill}).
Unlike classical billiards, the reflection law for symplectic billiards is not local.

In a recent paper \cite{baracco-bernardi}, using new ideas, it is proved that a totally integrable symplectic billiard is an ellipse. Moreover, in the centrally-symmetric case, the analog of Theorem \ref{thm:birkhoff-1/4} also holds true for symplectic billiards \cite{baracco-bernardi-nardi}.

The fourth billiard model is the \emph{outer length billiard} (sometimes called \emph{4th billiard}).
This model was introduced in \cite{albers-tabachnikov}, see also the recent preprint \cite{baracco-bernardi-corentin}.

\begin{figure}[h]
	\centering
	\includegraphics[width=0.65\linewidth]{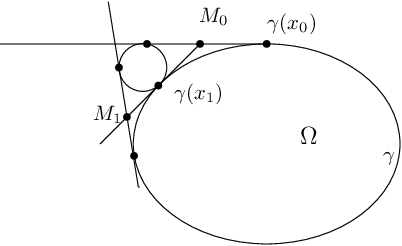}
	\caption{Outer length billiard rule, $B_{\Omega4}:M_0\mapsto M_1$.}
	\label{fig:fourth-rule}
\end{figure}

This billiard map, denoted by $ B_{\Omega4}$, acts in $\Omega^c$, the exterior of $\Omega$, according to the following rule, see Fig. \ref{fig:fourth-rule}.

In this figure, given a point $M_0\in\Omega^c$, consider two tangents through $M_0$ to $\gamma$ and the unique circle tangent to $\gamma$ at $\gamma(x_1)$ and to the tangent line at $\gamma(x_0)$. Then the point $M_1=B_{\Omega4}(M_0)$ is defined as the intersection point of the tangent line at $\gamma(x_1)$ and the unique line, which is tangent to the circle and $\gamma$. We will denote the lengths $$\lambda_0:=|M_0-\gamma(x_0)|,\quad \lambda_1:=|M_0-\gamma(x_1)|.$$
One can prove \cite{albers-tabachnikov} that $B_{\Omega4}$ is a positive twist map, and $S_{\Omega4}(x_0,x_1)=\lambda_0+\lambda_1$ is a generating function for the outer length billiard.

It was observed in \cite{baracco-bernardi-corentin}, using \cite{stachel}, that ellipses (and circles) are totally integrable outer length billiards. Moreover, the foliation of $\Omega^c$ by invariant curves is the foliation by confocal ellipses. However, the inverse rigidity result is not known in this case.

\section{Effective rigidity for Birkhoff billiards.}\label{sect:effective}
Rigidity theorems on circular and elliptic Birkhoff billiards allow effective versions. An effective rigidity of circular billiards was first found in \cite{bialy-effective}.

For the action functional $\mathcal L=\sum_n L(q_n,q_{n+1})$, we consider {\it locally maximizing} configurations, that is, those configurations which give the {\it local maximum} for the functional between any two endpoints.
We call such configurations \textit{m-configurations}, and the corresponding orbits on the phase cylinder $\mathbb{A}$ we call \textit{m-orbits}. We refer to \cite{bialy-effective,bialy-tsodikovich} for more details.

The set $\mathcal M$ of all m-orbits is a closed invariant subset, containing, in particular, all rotational invariant curves as well as Aubry-Mather sets. It was proved in \cite{bialy-effective} that the measure of the complement of $\mathcal{M}$ can be estimated from below in geometric terms (isoperimetric defect of the billiard curve).

Another way to see this kind of estimate uses the supporting function of a convex curve.
Given a smooth convex curve $\gamma$, its support function $h$ is a function $h:[0,2\pi]\to \mathbb R$, where $h(\psi)$ is the distance from the origin to the (oriented) tangent to $\gamma$, for which the right unit normal is $(\cos \psi,\sin \psi)$.
Write $\mu$ for the invariant measure on the phase cylinder of the billiard map.
\medskip

\begin{theorem}\label{thm:effectiveGeneralCurve}
	Let $\gamma$ be a planar strictly convex $C^2$-smooth curve with support function $h:[0,2\pi]\to\mathbb R$ (with respect to an arbitrary origin in the interior of $\gamma$).
	Let $\mathbb{A}$ denote the phase cylinder of the Birkhoff billiard map in $\gamma$, and let $\mathcal{M}\subseteq\mathbb{A}$ denote the set swept by m-orbits, and $\Delta=\mathbb{A}\setminus\mathcal{M}$.
	Then the following estimate holds true
	\[\mu(\Delta)\geq \pi^2k_{min} d^2(h,W),\]
	where $0<k_{min}$ is the minimal curvature of $\gamma$, $W$ denotes the subspace of $L^2[0,2\pi]$ spanned by the functions $\{1,\cos(\psi),\sin(\psi)\}$, and $d(\cdot,W)$ denotes the $L^2$-distance to that subspace.
	Moreover, this bound is sharp for circles.
\end{theorem}
\medskip

The method of the proof can also be adopted to give an effective version of the rigidity of integrable billiards in centrally-symmetric curves given by Theorem \ref{thm:birkhoff-1/4}.

Let us denote by $\mathcal{C}$ the set of centrally-symmetric, strictly convex, $C^2$-smooth curves, for which the billiard map has an invariant curve $\alpha$ of rotation number $1/4$ consisting of $4$-periodic orbits, as was considered in \cite{bialy-mironov}.
The class $\mathcal{C}$ is, in fact, rather big and can be fully characterized in terms of the support function $h$.
We will denote by $\mathcal{A}$ the domain bounded by $\alpha$ and the upper boundary of the cylinder, and by $\mathcal{M}$, the subset of $\mathbb{A}$ swept by m-orbits.
We give a lower bound for the measure of the invariant subset $\Delta_{\mathcal{A}}:=\mathcal{A}\setminus\mathcal{M}$, which is free from m-orbits.
\medskip

\begin{theorem}\label{thm:estimate-beta}
	Let $\gamma$ be a centrally-symmetric, strictly convex, $C^2$ smooth curve for which the billiard map has an invariant curve $\alpha$ of rotation number $1/4$ consisting of $4$-periodic orbits.
	Denote by $h:[0,2\pi]\to\mathbb R$ the support function of $\gamma$, with respect to the center of symmetry.
	Set $\Delta_{\mathcal{A}}=\mathcal{A}\setminus\mathcal{M}$.
	Then the following estimate for the measure holds true:
	$$
	\mu(\Delta_{\mathcal{A}})\geq \frac{25\pi^2}{32}k_{min}^3 d^2(h^2,U),
	$$
	where $0<k_{min}$ is the minimal curvature of $\gamma$, $U$ is the subspace of $L^2[0,\pi]$ spanned by $\{1,\cos(2\psi),\sin(2\psi)\}$, and $d(\cdot,U)$ is the $L^2$-distance from this subspace.
	Moreover, this bound is sharp for ellipses.
\end{theorem}
\medskip

\begin{remark}
	The bounds of Theorems \ref{thm:effectiveGeneralCurve}, \ref{thm:estimate-beta} can be considered as effective versions of the results on the Birkhoff--Poritsky conjecture \cite{B0, bialy-mironov} by estimating the closeness of $\gamma$ to the ``closest" circle (for arbitrary curves) or ellipse (for centrally-symmetric curves) in terms of the measure of the set $\Delta$.

\end{remark}

\begin{figure}[h]
	\centering
	\includegraphics[width=0.55\textwidth]{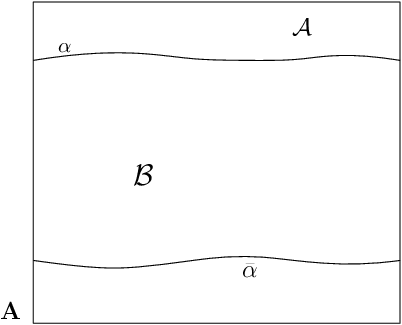}
	\caption{Region $\mathcal B$ between the invariant curves $\alpha,\bar\alpha$; $ \Delta_\mathcal B=\mathcal B\setminus \mathcal M, \Delta_\mathcal A=\mathcal A\setminus \mathcal M$}
	\label{maintheorem}
\end{figure}

It was shown (see \cite{effective-away}) that an effective sharp estimate can be obtained also for a portion of the set $\mathcal M$ lying in the domain $\mathcal B$ bounded by two invariant curves $\alpha$ and $\bar\alpha$ of 4-periodic orbits. Here, $\bar\alpha$ is the invariant curve of rotation number $\frac{3}{4}$ which corresponds to $\alpha$. This curve consists of the same billiard trajectories as $\alpha$, but with the reversed orientation of the lines, see Fig. \ref{maintheorem}. 
\medskip

\begin{theorem}\label{main1}
	Let $\gamma\subset\mathbb{R}^2$ be a $C^2$-smooth, centrally-symmetric, convex closed curve of positive curvature. We will assume that the billiard map corresponding to $\gamma$ has a rotational (i.e., winding once around the cylinder and simple) invariant curve $\alpha\subset \mathbb{A}$ consisting of $4$-periodic orbits.
	Let $\bar\alpha$ be the corresponding invariant curve of rotation number $\frac{3}{4}$.
	Let $\mathcal B\subset \mathbb A$ be the domain between the curves $\alpha$ and $\bar\alpha$.
	Then the following estimate holds:
	$$
	\frac{3k_{min}}{16}(P^2-4\pi A)\leq \mu(\Delta_{\mathcal{B}}),
	$$
	where $P, A$ denote the perimeter and the area of $\gamma$, $\Delta_{\mathcal{B}}=\mathcal{B} \setminus \mathcal{M}$, and $k_{min}>0$ is the minimal curvature of $\gamma$.
\end{theorem}
\medskip

Notice that the region $\mathcal B$ lies apart from the boundary of the phase cylinder. Nevertheless we have:
\begin{corollary}\label{circle}
	The set $\mathcal M$ of locally maximizing orbits contains the whole region $\mathcal B$ if and only if $\gamma$ is a circle.
\end{corollary}

\section{Billiards with finite group of rotational symmetry.}\label{sect:symmetry}
In this section we describe the results of \cite{bialy-imrn} on the rigidity by one invariant curve for billiard tables with a finite rotational symmetry group. The motivation comes from Theorem \ref{thm:birkhoff-1/4}. Namely, we consider Birkhoff billiard tables with rotational symmetry of order $k\geq 3$. We show that among such tables, the only ones that have a rotational invariant curve that consists of $k$-periodic orbits are circles (an exact formulation will appear below). Note that this claim does not assume foliation by invariant curves, but it assumes the existence of a specific invariant curve, the one consisting of $k$-periodic orbits. A similar result can be proved for the outer length billiard.

The main results are the following:
\medskip

\begin{theorem}\label{theorem euclideanK}
	Let $\gamma$ be a $C^2$-smooth, planar, strictly convex curve which is invariant under $R_{\frac{2\pi}{k}}$, with $k\geq 3$.
	If the Birkhoff billiard map inside $\gamma$ has a rotational invariant curve of $k$-periodic orbits, then $\gamma$ is a circle.
\end{theorem}
\medskip

Similarly, we state for the outer length billiard (it is not contained in \cite{bialy-imrn}, but the proof is completely analogous):
\medskip

\begin{theorem}\label{theorem outerlength}
	Let $\gamma$ be a $C^2$-smooth, planar, strictly convex curve which is invariant under $R_{\frac{2\pi}{k}}$, with $k\geq 3$.
	If the outer length billiard map outside $\gamma$ has a rotational invariant curve of $k$-periodic orbits, then $\gamma$ is a circle.
\end{theorem}
\medskip

Moreover, we can state similar results for outer and symplectic billiards. For Minkowski billiards (see, e.g., \cite{GUTKIN2002277}) we will obtain a slightly different result which we discuss in Section \ref{sect:magnetic}. Since these billiards are equivariant with respect to all affine transformations, we assume in these settings that the billiard table is invariant under a (finite) order $k\geq 3$ element of $\textrm{GL}(2,\mathbb R)$. Thus, for outer and symplectic billiards, we show that if for such a domain there exists a rotational invariant curve of $k$-periodic orbits, then the domain must be an ellipse. It is an elementary fact in linear algebra that an order $k\ge 3$ element in $\textrm{GL}(2,\mathbb R)$ is conjugate to a rotation by angle $\frac{2\pi m}{k}$, where $m$ and $k$ are co-prime. For such $k,m$, any set that is invariant under a rotation by $\frac{2\pi m}{k}$ is also invariant under a rotation by $\frac{2\pi}{k}$. For this reason, we can always assume that the order $k$ linear map is conjugate to a rotation by angle $\frac{2\pi}{k}$. 
Here and below, we write $R_\theta$ for the rotation of $\mathbb R^2$ by angle $\theta$.
\medskip

\begin{theorem}\label{theorem outer}
	Let $\gamma$ be a $C^2$-smooth, planar, strictly convex curve which is invariant under a linear map of order $k\ge 3$.
	If the outer billiard map in the exterior of $\gamma$ has a rotational invariant curve of $k$-periodic orbits, then $\gamma$ is an ellipse.
\end{theorem}
\medskip

Here the restriction $k\ge 3$ is natural: there are no $2$-periodic orbits for the outer billiard map.

Similarly, we state for the symplectic billiard:
\medskip

\begin{theorem}\label{theorem symplectic}
	Let $\gamma$ be a $C^2$-smooth, planar, strictly convex curve which is invariant under a linear map of order $k\ge3$.
	If the symplectic billiard map inside $\gamma$ has a rotational invariant curve of $k$-periodic orbits, then $\gamma$ is an ellipse.
\end{theorem}
In this theorem, the assumption that $k\ge 3$ is necessary: by the definition of symplectic billiards, every curve $\gamma$ admits an invariant curve of two periodic orbits: these are the orbits between the pairs of points of $\gamma$ that have parallel tangents.

Now we turn to an important application of Theorem \ref{theorem outer} to the geometry of Radon planes. A normed plane $( \mathbb {R}^2, ||\cdot||)$ is called a Radon plane if the Birkhoff orthogonality relation is symmetric (see \cite{martini2006antinorms}). The unit circle of a Radon norm is called a Radon curve. An equivalent definition in terms of the outer billiard map (see \cite {bialy2022self}) reads: {\it Let $\gamma$ be the unit circle of the $C^2$-smooth norm, then the norm is Radon if and only if the outer billiard map corresponding to $\gamma$ has a rotational invariant curve of 4-periodic orbits.} Using the last definition, we get an immediate corollary of Theorem \ref{theorem outer}:
\begin{corollary}
	If a Radon norm on $\mathbb {R}^2$ is invariant under a linear map of order $4$, then the norm is necessarily Euclidean.
\end{corollary}
\medskip

\begin{figure}[h]
	\centering
	\includegraphics[width=1\textwidth]{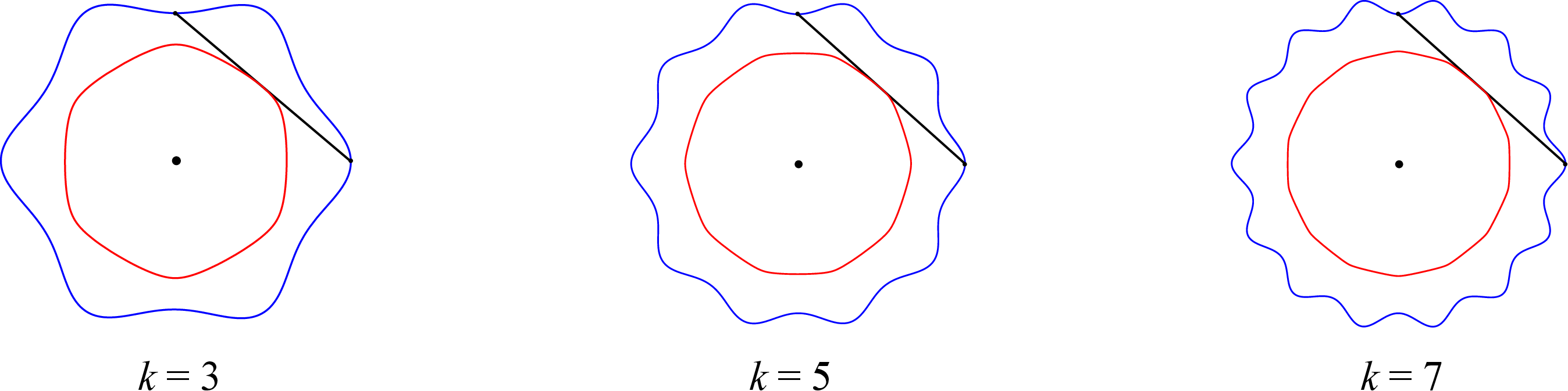}
	\caption{Radon curves (red) and invariant curve of 4-periodic orbits (blue) having symmetries of order $6,10,14$ respectively}
	\label{fig:radonExamples}
\end{figure}

Remarkably, examples of non-circular analytic Radon curves having symmetries of order $2k$, with odd $k\geq3$, were constructed using Lame differential equation in a recent paper, see  \cite{bialy2022self}, and Fig. \ref{fig:radonExamples}.

\section{Mather $\beta$-function and integrability.}\label{sect:beta}
Given a positive twist exact symplectic  map of the cylinder, the Mather $\beta$-function assigns to a rotation number $\rho$ the average action of a minimal trajectory with this rotation number. This function encodes the action spectrum of the twist map.

\medskip
In particular, one has for the billiard models the following results.
\medskip

\begin{example}
	\begin{enumerate} 
		\item[]
		\item For Birkhoff billiard, $\beta(\frac{m}{n})=-$ {\it the maximal perimeter of the periodic billiard configuration divided by the number of vertices.}
		\item For outer billiard, $\beta(\frac{m}{n})=$ {\it the area covered by the trajectory divided by the number of vertices.}
		\item For symplectic billiard, $\beta(\frac{m}{n})=-$ {\it the area covered by the trajectory divided by the number of vertices.}
		\item For outer length billiard, $\beta(\frac{m}{n})=$ {\it the minimal perimeter of the periodic billiard configuration divided by the number of vertices.}
	\end{enumerate}
\end{example}

Notice that one needs to introduce minus sign for the generating functions of Birkhoff and symplectic billiards in order to have a positive twist map. Hence, the minus sign appears for the corresponding $\beta$-function for Birkhoff and symplectic billiards in the examples.

J. Mather proved that the function $\beta$ is continuous and strictly convex. Moreover, $\beta$ is differentiable at any $\rho\in\mathbb R\setminus \mathbb Q$. The function $\beta$ is differentiable at  $\rho\in\mathbb Q$ if and only if there is a rotational invariant curve consisting of periodic orbits of rotation number $\rho$.

It is very natural to use the Mather $\beta$-function for inverse geometric problems of billiards. We refer to   \cite{KSH} and to a very detailed survey \cite{Kaloshin-lecturenotes} on various results and techniques.

For example, Theorems \ref{theorem euclideanK}, \ref{theorem outerlength}, \ref{theorem outer} and \ref{theorem symplectic} can be stated uniformly as follows:
\medskip

\begin{theorem}
	Let $\gamma\subset\mathbb R^2$ be a convex $C^2$-smooth closed curve of positive curvature. Assume $\gamma$ enjoys symmetry of order $k$ as formulated in Theorems \ref{theorem euclideanK}, \ref{theorem outerlength}, \ref{theorem outer} and \ref{theorem symplectic}. If the Mather $\beta$-function is differentiable at a rational number $\frac rk$, for co-prime $r,k$, then $\gamma$ is a circle in the cases of Birkhoff and outer length billiards, but $\gamma$ is an ellipse in the cases of outer or symplectic billiards.
\end{theorem}
\medskip

Now, we would like to stress the following corollaries of the theorems on total integrability discussed above. The idea to use the Mather differentiability theorem in this context belongs to K. F. Siburg \cite{siburg} (and before that, by V. Bangert for Riemannian tori).

The following (see \cite{bialy-entropy}) can be deduced from Theorem \ref{thm:birkhoff-1/4}:
\medskip

\begin{theorem}
	Let $\gamma=\partial \Omega$ be a simple, closed, centrally-symmetric curve of positive curvature. If Mather's function $\beta$ is differentiable on $ (0, 1/4]$, then $\gamma$ is an ellipse.
\end{theorem}
\medskip

\begin{corollary}
	Let $\Omega_1, \Omega_2$ be two strictly convex $C^2$-smooth centrally-symmetric planar domains such that $\Omega_1$ is an ellipse.
	Suppose that the Mather $\beta$-functions $\beta_1,\beta_2$ of the domains satisfy
	
	$${\beta_1(\rho)=\beta_2(\rho),\  \forall \rho\in \left(0,\frac{1}{4}\right].}$$
	
	Then $\Omega_2$ is an ellipse isometric to $\Omega_1$.
\end{corollary}
\medskip

In order to prove the corollary, namely that the $\Omega_1,\Omega_2$ are isometric, one shows that it is possible to recover an ellipse by the two values $\beta'(0)$ and $\beta(\frac 14)$.

The following question naturally arises: How many values of the Mather $\beta$-function determine the ellipse in the class of ellipses. Can one recover an ellipse by any two values of the $\beta$-function $\beta(\rho_1), \beta(\rho_2)$ for the rotation numbers $ \rho_{1,2}\in(0,\frac{1}{2}]$?

One can prove (see \cite{bialy-entropy}) that an elliptic Birkhoff billiard can be determined by two values of $\beta(\rho_1),\beta(\rho_2)$, where $\rho_1=\frac{1}{2}$ and $\rho_2=\frac{m}{n}$ is any rational in $(0,\frac{1}{2})$. Notice that though the Mather $\beta$-function for an elliptic Birkhoff billiard can be explicitly found via elliptic functions (see \cite{bialy-entropy}), this does not help much in the proof of this kind of results and therefore new geometric arguments are needed.
\bigskip

For outer billiards, we can conclude in a similar way from Theorem \ref{thm:total} that an elliptic outer billiard table can be recovered by the Mather $\beta$-function (i.e., by the circumscribed marked Area spectrum) for $\rho\in(0;1/2)$:
\medskip

\begin{theorem}
	Let $\gamma=\partial\Omega$ be a convex curve of positive curvature.
	If the Mather $\beta$-function of the outer billiard is differentiable on $(0,1/2)$ then $\gamma$ is an ellipse.
\end{theorem}
\medskip

It is an open question whether for a centrally-symmetric outer billiard table one can relax differentiability of the Mather $\beta$-function on $(0,1/2)$, to require the differentiability on $(0,1/4]$ or on $[1/4, 1/2)$.

However, for symplectic billiards we can conclude from the results of \cite{baracco-bernardi} and \cite{baracco-bernardi-nardi}:

\begin{corollary}
	It follows from \cite{baracco-bernardi} and \cite{baracco-bernardi-nardi} that symplectic billiard in ellipse can be recovered by the Mather $\beta-$function (i.e., by the inscribed marked Area spectrum) for $\rho\in(0;1/2)$, or $\rho\in(0;1/4)$ in the centrally-symmetric case.
\end{corollary}
\medskip

We turn now to state an isoperimetric inequality for the Mather $\beta-$ function, which is related to the approximation of convex curve by inscribed/circumscribed extremal polygons. This inequality was found very recently in \cite{baranzini-bialy-sorrentino}. It can be proved for three (out of four) billiard models discussed in Section \ref{sect:other}, namely for Birkhoff, symplectic and outer length billiard models, but surprisingly not for the outer billiard model, which is exceptional.
\medskip

\begin{theorem}\label{thm:isoperimetric}
	Let $\Omega$ be a convex billiard table, and let $\partial \Omega$ be a $C^2$-smooth curve of positive curvature. Then for the $\beta$-function of Birkhoff, symplectic, and outer length billiards, the following inequality holds:
	$$
	\beta_\Omega(\rho)\leq\beta_D(\rho).
	$$
	Here $D\subset\mathbb R^2$ is a disc of the same boundary length as $\Omega$ in the cases of Birkhoff and outer length billiards, but of the same area as $\Omega$ in the case of symplectic billiard.
\end{theorem}
\medskip

Surprisingly, Theorem \ref{thm:isoperimetric} does not necessarily hold for outer billiards. It is shown in \cite{baranzini-bialy-sorrentino} using Santalo and Blaschke inequalities, that for Radon curves the sign in the inequality of Theorem \ref{thm:isoperimetric} is always reversed! The sign is also reversed for outer billiards having invariant curve of 3-periodic orbits.

Interestingly, the equality case in Theorem \ref{thm:isoperimetric} is related to the so called Gutkin billiard tables \cite{gutkin}. These are Birkhoff billiards which have an invariant curve of constant reflection angle $\pi\rho$. Moreover, non-circular Gutkin billiard table exist only when $\rho$ is a solution of $\tan(k\pi\rho)=k\tan(\pi\rho)$ for some $k\geq 3$.
\medskip

\begin{theorem}
	\begin{itemize}
		\item[] 
		\item [(1)]{If for some rotation number  $\rho\in \mathbb Q$ the equality in Theorem \ref{thm:isoperimetric} holds, then $\Omega$ is a circle in the cases of Birkhoff and outer length billiards, and $\Omega$ is an ellipse for symplectic billiards.}
		\item [(2)]{If for $\rho\in \mathbb R\setminus \mathbb Q$ the equality in Theorem \ref{thm:isoperimetric} holds for a Birkhoff billiard, then $\Omega$ is a Gutkin billiard table.}
		\item [(3)]{If for $\rho\in \mathbb R\setminus \mathbb Q$ the equality in Theorem \ref{thm:isoperimetric} holds for the cases of outer length or symplectic billiards, then $\Omega$ is a circle in the case of outer length billiard, and $\Omega$ is an ellipse in the case of symplectic billiard.}
	\end{itemize}
\end{theorem}
\medskip

\section{Magnetic and Minkowski billiards.}\label{sect:magnetic}
In this section we discuss separately some results on two very interesting billiard models, namely on magnetic and Minkowski billiards.

\begin{figure}[h]
	\centering
	\includegraphics[width=0.9\textwidth]{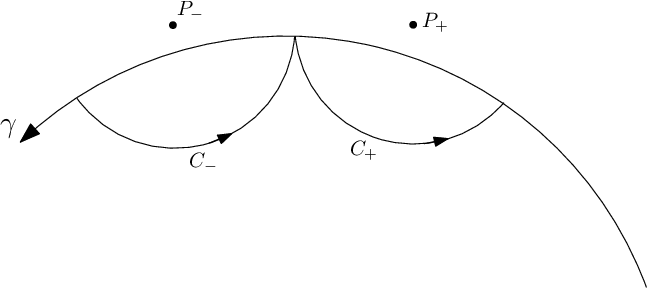}
	\caption{Reflection rule for Larmor circles and their centers}
	\label{fig:mag}
\end{figure}

Magnetic billiards were introduced and studied in \cite{berry-robnik} and since then attracted many researchers. In this model, the billiard ball moves under the influence of a magnetic field along Larmor circles (for the constant magnetic field of magnitude $B$, the Larmor circles are of radius $1/B$), see Fig. \ref{fig:mag}.

The only known example of an integrable magnetic billiard is a circular billiard in a constant magnetic field. We believe that for non-constant magnetic fields there exist new integrable examples. (For comparison, there are many new integrable magnetic geodesic flows on the 2-torus, \cite{magnetic-geodesic}.)

We have found that the algebraic approach can be extended to the case of magnetic billiards in a constant magnetic field. This extension is based on the interplay of differential and algebra-geometric properties of the equidistant curves of the boundary of the billiard domain. We have implemented it in \cite{BM3,BM5} for the plane magnetic billiards and then in \cite{BM4} for magnetic billiards on the constant curvature surfaces.

Moreover, the approach of total integrability works well for convex magnetic billiard in a {\it weak} constant magnetic field \cite{bialy-magnetic-hopf}. This is done using the magnetic mirror equation (generalizing the approach suggested in \cite{W}) and the magnetic Santalo formula.

Consider the magnetic billiard inside a convex curve $\gamma$, where the magnitude of the magnetic field is assumed to be constant and relatively small with respect to the geodesic curvature of the boundary,
$0<B < \min_{x\in \gamma}k(x).$
\medskip

\begin{theorem}\label{main}
	If the magnetic billiard in $\gamma$, which is subject to a weak constant magnetic field, is totally integrable, then $\gamma$ must be a circle.
\end{theorem}
\medskip

It is not clear how the assumption of a weak magnetic field in this theorem can be relaxed.
\medskip

The model of Minkowski and, more generally, Finsler billiards, was introduced and studied by Gutkin and Tabachnikov in \cite{GUTKIN2002277}. It plays an important role in a recent study by Artstein-Avidan, Karasev, and Ostrover \cite{Artstein_Avidan_2014} on Mahler's conjecture. Consider a smooth convex curve $\gamma\subset\mathbb R^2$, and a Minkowski (not necessarily symmetric) norm $N$. The Minkowski billiard in $\gamma$ with respect to the norm $N$ is defined as a twist map that corresponds to the generating function \[L(t_0,t_1)=N(\gamma(t_1)-\gamma(t_0))\,,\] and the corresponding variational principle (\ref{functional}), where $\gamma$ is parametrized by $N$-arc length parameter (meaning, $N(\dot{\gamma}(t))=1$). Thus, the Birkhoff billiard is a special case of the Minkowski billiard when $N$ is chosen to be the Euclidean norm.

An important special case of Minkowski billiard occurs when $\gamma$ is the unit ball of the norm $N$.

The following is an interesting question: what are integrable Minkowski billiards? In particular, is it true that a totally integrable Minkowski billiard, corresponding to the norm $N$, in a unit ball of $N$ is necessarily Euclidean?

This question was approached in \cite{bialy-imrn} for billiards with a finite symmetry group. In particular, the following is proved in \cite{bialy-imrn}:
\medskip

\begin{theorem}\label{cor notTotalIntegrableMinkowski}
	Let $3\le k\in\mathbb N$, and let $\gamma$ be a $C^2$-smooth, planar strictly convex curve which is invariant under a linear map of order $k$.
	Suppose that the Minkowski billiard map in $K$, where $K$ is the interior of $\gamma$, with the norm induced by $K$, has a rotational invariant curve of rotation number $\frac{r}{k}$, with $r$ co-prime to $k$.
	If this invariant curve has a neighborhood $U$ in which the Minkowski billiard map is totally integrable, then $K$ is also invariant under a linear map of order $ak$, where $a$ depends on the remainder of $k$ modulo $4$:
	\[\begin{cases}
		a = 1\,, \textrm{if } k\equiv 2 \pmod{4} \,, \\
		a = 2\,, \textrm{if } k\equiv 0 \pmod{4} \,, \\
		a = 4\,, \textrm{if } k\equiv 1 \pmod{2} \,. \\
	\end{cases}\]
\end{theorem}
\medskip

Moreover, Theorem \ref{cor notTotalIntegrableMinkowski} can be ``iterated'' to obtain the following rigidity result for Minkowski billiards, under the symmetry assumption:
\medskip

\begin{theorem} \label{thm BirkhoffForMinkowski}
	Let $\gamma$ be a $C^2$-smooth, planar, strictly convex curve which is invariant under a linear map of order $k\ge 3$.
	Consider the Minkowski billiard system in $K$, the interior of $\gamma$, with the norm induced by $K$.
	If the Minkowski billiard map of $K$ is totally integrable on the entire phase cylinder, then $\gamma$ is a (Euclidean) ellipse.
\end{theorem}
\medskip

\begin{example}
	Consider the $L^p$-balls in $\mathbb R^2$ for $p> 2$.
	It is well known that the only linear isometries of the $L^p$-norm, for $p\neq 2$ are generalized permutations (that is, products of permutation matrices with diagonal matrices with diagonal entries $\pm{1}$).
	There exists a generalized permutation of order four, for example, rotation by $\frac{\pi}{2}$, but in $\mathbb R^2$ there are no generalized permutations of order eight.
	By Theorem \ref{cor notTotalIntegrableMinkowski}, the Minkowski billiard map inside an $L^p$-ball is not totally integrable.
\end{example}
\medskip

\section{Wire billiards and billiards in cones.}\label{sect:wire}

In this section we first consider the wire billiards introduced and studied in \cite{wire,DM}, and then billiards inside smooth cones ${\mathbb R}^n$ \cite{MY,MY1}.

We turn now to describe the examples of integrable wire billiards, admitting first integrals which are polynomial in the components of the velocity vector.

Let $\gamma\subset {\mathbb R}^n$ be a smooth curve. The wire billiard is defined as follows. Let $\gamma(s),\gamma(s_1)$ be two points on the curve. The chord $[\gamma(s_1)\gamma(s)]$ is reflected to a chord $[\gamma(s)\gamma(s_2)]$ if the angles between the chords and the tangent vector $\dot \gamma(s)$ are equal. For a small perturbation of the plane curve, the point $\gamma(s_2)$ is defined uniquely, but not necessarily in general. In the general case, the wire billiard map can be multi-valued.

In \cite{wire}, an example of integrable wire billiards is given.
Consider the wire curve
\begin{equation}\label{eq:2}
	\gamma(s)=e^{As}\gamma_0,\ A\in so(n),\ \gamma_0\in{\mathbb R}^n.
\end{equation}
One can check that for any chord $[\gamma(x),\gamma(y)]$, the angle between the chord and the vectors $\dot\gamma(x)$ and $\dot\gamma(y)$ are equal. Therefore, this angle is a first integral, i.e., it is preserved under the dynamics. So, this gives an example of totally integrable wire billiard.
It can be checked \cite{DM} that conservation of the angle leads to a polynomial integral of degree one in the velocity. Namely, let $v=(v_1,\dots,v_n)$ be the velocity vector of the particle, $|v|=1,x=(x_1,\dots,x_n)$ are coordinates of the particle. Then the function $F$ is a polynomial integral of the wire billiard flow,
$$
F(x,v) =\langle Ax,v \rangle=\sum_{i<j}a_{ij} (v_jx_i - v_ix_j),
$$
where $a_{ij}$ are the components of the matrix $A$.

In the case of ${\mathbb R}^3$, the curve $\gamma(s)$ in (\ref{eq:2}) is a circle in a 2-plane. It would be interesting to construct a closed, non-planar curve in ${\mathbb R}^3$ such that the corresponding wire billiard is  integrable. Similarly to (\ref{eq:2}), integrable round spirals can be constructed, but they are not closed.

In the case of ${\mathbb R}^{n}, n>3$, there exist non-planar closed curves $\gamma$ given by (\ref{eq:2}).

\medskip
\begin{example}
	The wire billiard corresponding to a closed curve $\gamma$ in ${\mathbb R}^4={\mathbb C}^2$ given by the formula
	$$
	\gamma(t)=(ae^{ikt},be^{imt}), a,b>0, t\in[0,2\pi],k,m\in{\mathbb N}
	$$
	is totally integrable. The curve is a toric knot in $S^3\subset{\mathbb R}^4.$ Moreover, taking $a=1, b=\varepsilon,k=1$, we get a deformation of the planar circle to totally integrable wire billiards.
\end{example}
\medskip

We turn now to state that the Birkhoff billiard inside a cone over $C^3$-smooth, strictly convex closed hypersurface embedded in ${\mathbb R}^{n-1}$ with non-degenerate second fundamental form admits first integrals whose values uniquely determine all billiard trajectories.

It is well known \cite{tabachnikov1995billiards} that the billiard trajectory in the planar angle $\alpha$ has at most
$\left\lceil \frac{\pi}{\alpha} \right\rceil$ reflections, where $\left\lceil x \right\rceil$ is the smallest integer greater than or equal to $x$.

Y. G. Sinai \cite{Sin} proved that every trajectory in a polyhedral cone also has a finite number of reflections. Moreover, the number of reflections is bounded by a universal constant depending only on the cone.

An analog of Sinai's theorem in the smooth case was proved in \cite{MY}.
\medskip

\begin{theorem}\label{cone-bounded}
	Let $K$ be the cone defined by
	$$
	K = \{ tp \mid p \in \Gamma, t > 0 \} \subset \mathbb{R}^n,
	$$
	where $\Gamma \subset \{ x \in \mathbb{R}^n \mid x^n = 1 \}$ is a $C^3$-smooth, strictly convex, closed submanifold with a positive definite second fundamental form at every point.
	Then any billiard trajectory in $K$ has a finite number of reflections.
\end{theorem}
\medskip

Notice however, that for smooth cones the number of reflections cannot be uniformly bounded.

It is an interesting question, for what submanifolds $\Gamma\subset{\mathbb R}^{n-1}$ any billiard trajectory inside the cone over $\Gamma$ has a finite number of reflections?

In Theorem \ref{cone-bounded}, the assumption of $C^3$-smoothness of $\Gamma$ is essential. Namely, we have (see \cite{MY1}):
\medskip

\begin{theorem}
	There exist $C^2$-smooth convex cones with billiard trajectories having infinitely many reflections in finite time.
\end{theorem}
\medskip

M. Berger \cite{Ber} proved that if the Birkhoff billiard for $C^2$-smooth hypersurface $\Sigma \subset \mathbb{R}^n,$ $n \geq 3$  with non-degenerate second fundamental form admits a $C^2$ caustic with non-degenerate second fundamental form, then $\Sigma$ must be a piece of a quadric and the caustic a piece of a confocal quadric. In the case of a cone, whose second fundamental form is positive semi-definite, the analog of Berger's theorem is not valid.
Namely, for any cone, there always exists a family of caustics that are spheres \cite{MY}.
\begin{figure}[htbp]
	\centering
	\includegraphics[scale=0.09]{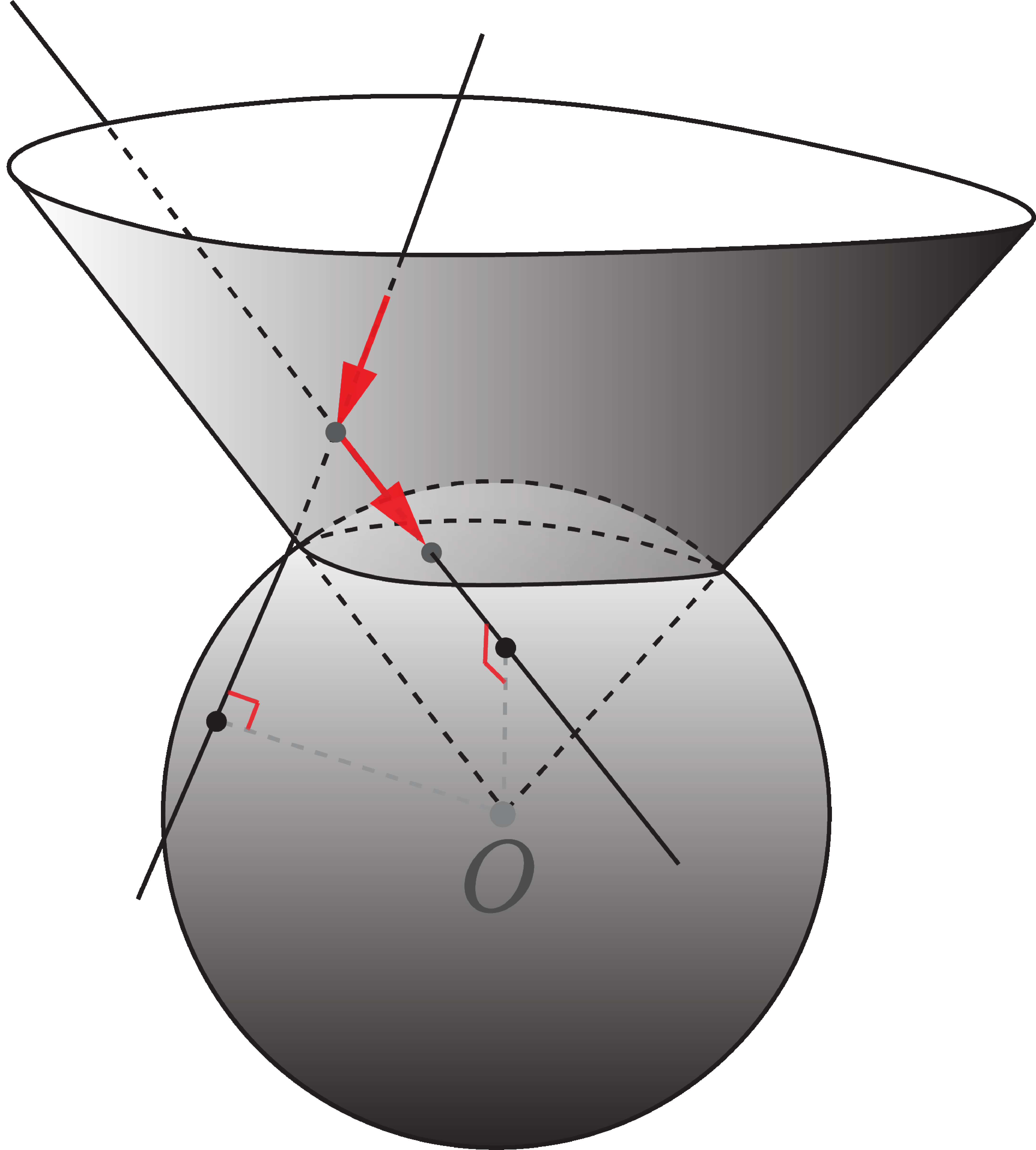}
	\caption{The sphere as a caustic of the billiard inside a cone.}
	\label{fig:caustic}
\end{figure}
\medskip

\begin{theorem}
	\begin{enumerate}
		\item[]
		\item [(1)]The Birkhoff billiard inside $K \subset \mathbb{R}^n$ admits the first integral
		$$
		I = \sum_{1 \leq i<j\leq n} m_{i,j}^2,
		$$
		where $m_{i,j} := x^i v^j - x^j v^i$ for $i < j$, $i, j = 1, \ldots, n$, and $v = (v^1, \ldots, v^n)$ is the velocity vector.
		\item [(2)]The spheres centered at the vertex $O \in \mathbb{R}^n$ of $K$ are caustics of the billiard inside $K$ (see Fig. \ref{fig:caustic}).
	\end{enumerate}
\end{theorem}
\medskip

The phase space $\Psi$ of the Birkhoff billiard inside $K$ is an open subset of $T\mathbb S^{n-1}$ consisting of oriented lines in $\mathbb R^n$ that intersect $K$ transversally. Let ${B_K}:\Psi\rightarrow\Psi$ be the billiard map.
\medskip

\begin{theorem}
	There are continuous (smooth almost everywhere) first integrals
	$I_1(x)$, $\ldots$, $I_{2n-2}(x)$, $I_{2n-1}(x)=I$ on $\Psi$
	invariant under $B_K$,
	$$
	I_j(x)= I_j(B_K(x)), \quad  x \in  \Psi,
	\quad j=1, \ldots, 2n-1.
	$$
	Each point in the image of the map
	$$
	\mathcal{I}= (I_1,\ldots, I_{2n-1}): \Psi \to \mathbb{R}^{2n-1}
	$$
	determines a unique billiard trajectory in $\Psi$.
\end{theorem}
\medskip

There is a smooth submanifold ${\sigma} \subset \Psi$ such that $I_j$ are smooth on $\Psi\setminus \sigma$. There exist also {\it smooth} first integrals $I^s_1(x)$, $\ldots$, $I^s_{2n-1}(x)$ on $\Psi$. The values of these integrals uniquely determine billiard trajectories in $\Psi \setminus \sigma$ and vanish identically on $\sigma$.

\section{Open questions.}\label{sect:questions}
In this section we formulate natural open questions related to the results discussed in previous sections.

\subsection{Birkhoff billiards.}
\begin{enumerate}
	\item [(1)] Is it possible to remove the central-symmetry assumption in Theorem \ref{thm:birkhoff-1/4}?
	\item [(2)] Can one claim the result of Theorem \ref{thm:birkhoff-1/4} for a smaller neighborhood of the boundary? Maybe for billiard tables of higher symmetry?
	\item [(3)] Suppose for the Birkhoff billiard in $\gamma$ there exist a sequence $\alpha_n\subset\mathbb A$ of rotational invariant curves with rotation numbers tending to $\frac 12$. Is it true that in this case $\gamma$ is an ellipse? If all invariant curves $\alpha_n$ correspond to convex caustics, then this is true \cite{BA}.
\end{enumerate}
\subsection{Outer billiards.}
\begin{enumerate}
	\item [(1)]  Are there new integrable examples of outer billiards?
	\item [(2)] Can one establish rigidity for total integrability on a part of the phase space?
	\item [(3)] In the centrally-symmetric case is it possible to state analog to Theorem \ref{thm:birkhoff-1/4}?
\end{enumerate}
\subsection{Outer length billiards.}
\begin{enumerate}
	\item [(1)]  For outer length billiards, ellipses are totally integrable. Are there other totally integrable billiards?
\end{enumerate}
\subsection{Magnetic billiards.}
\begin{enumerate}
	\item [(1)]  Is it possible to state an analog of Theorem \ref{main} for strong constant magnetic fields? The algebraic approach hints that this is indeed the case.
	\item [(2)] It would be interesting to find new integrable examples of magnetic billiards when the magnetic field is allowed to be variable.
\end{enumerate}
\subsection{Minkowski billiards.}
\begin{enumerate}
	\item [(1)] Prove or disprove that a Minkowski billiard corresponding to the norm $N$ inside the unit ball of $N$ is totally integrable if and only if $N$ is a Euclidean norm (communicated by Y.Ostrover).
	\item [(2)] Find new integrable examples of Minkowski billiards.
\end{enumerate}
\subsection{Mather $\beta$-function.}
\begin{enumerate}
	\item [(1)] Can one prove the reversed isoperimetric inequality for the Mather $\beta$-function of outer billiards, assuming the existence of a rotational invariant curve. As we mentioned above, this is the case for rotation numbers $\frac 13, \frac 14$.
	\item [(2)] Can one reconstruct an elliptic Birkhoff billiard in the class of ellipses by two values of $\beta$-function?
\end{enumerate}
\subsection{Wire and cone billiards.}
\begin{enumerate}
	\item [(1)] Are there integrable examples of wire billiards in $\mathbb R^3$? In \cite{wire}, integrable examples of closed wires are found in $\mathbb R^{n}, n>3$?
	\item [(2)] Can one deform in $\mathbb R^4$ an ellipse $E\subset\mathbb R^2\subset \mathbb R^4$ so that the wire billiard remains integrable?
	\item [(3)] For which submanifolds $\Gamma\subset{\mathbb R}^{n-1}$ any trajectory inside the cone $K$ over $\Gamma$ has finite number of reflections? Can one relax the condition of positive definiteness of the second fundamental form in Theorem \ref{cone-bounded}? For example, can the second fundamental form of $\Gamma$ be positive semi-definite, or indefinite? Can $\Gamma$ be of a topological type different from sphere, for example, of a torus?
\end{enumerate}
\bibliography{bibliography}
\bibliographystyle{abbrv}

Misha Bialy, 

School of Mathematical Sciences, Tel Aviv University, Israel;
	
e-mail: bialy@tauex.tau.ac.il 

${ \ }$

Andrey E. Mironov,

Sobolev Institute of Mathematics and Novosibirsk State University, Russia;

e-mail: mironov@math.nsc.ru

\end{document}